\input amstex
 \documentstyle{amsppt}
 \magnification=\magstep1
 \vsize=24.2true cm
 \hsize=15.3true cm
 \nopagenumbers\topskip=1truecm
 \headline={\tenrm\hfil\folio\hfil}

 \TagsOnRight

\hyphenation{auto-mor-phism auto-mor-phisms co-homo-log-i-cal co-homo-logy
co-homo-logous dual-izing pre-dual-izing geo-metric geo-metries geo-metry
half-space homeo-mor-phic homeo-mor-phism homeo-mor-phisms homo-log-i-cal
homo-logy homo-logous homo-mor-phism homo-mor-phisms hyper-plane hyper-planes
hyper-sur-face hyper-sur-faces idem-potent iso-mor-phism iso-mor-phisms
multi-plic-a-tion nil-potent poly-nomial priori rami-fication sin-gu-lar-ities
sub-vari-eties sub-vari-ety trans-form-a-tion trans-form-a-tions Castel-nuovo
Enri-ques Lo-ba-chev-sky Theo-rem Za-ni-chelli in-vo-lu-tion Na-ra-sim-han Bohr-Som-mer-feld}

 % reduced subscheme
 % minimal surface
\define\rest#1{_{\textstyle{|}#1}} % restriction of map to subset

\define\Span#1{\left<#1\right>} % <x> span or hull of x

 % Volume
 % roman differential dx, dy
 % identity map

% Bourbaki letters
 % complex numbers
\define\R{\Bbb R} % real numbers
\define\Z{\Bbb Z} % integers

 % projective space
 % projective plane
 % dual projective plane

% Script letters
\define\sA{{\Cal A}} % sheaf of algebras A
\define\sB{{\Cal B}} % conductor ideal
 % sheaf G
 % sheaf J
 % moduli space Mg
 % structure sheaf

% short Greeks

\define\om{\omega}

\define\Om{\Omega}

% \mathops
 % adjoint bundle

 % cardinality
 % codimension
 % cokernel of a map

 % image of a map
 % rank of a map

 % modulo
 % base locus of linear system

 % Grassmann variety
 % Hilbert scheme
 % Hom group
 % End group

 % Picard scheme
 % singular locus
 % Stabiliser group
 %second exterior power of a module

% Lie groups

 % unitary group

 \document

  \topmatter
  \title Hamiltonian dynamics on the moduli space of half weighted Bohr- Sommerfeld Lagranigian subcycles of a fixed volume \endtitle
  \author Nik. Tyurin                    \endauthor

   \address MPI, Bonn
   \endaddress
  \email
    tyurin\@tyurin.mccme.ru  jtyurin\@mpim-bonn.mpg.de
   \endemail

\abstract Connecting ideas of [2], [8] and [1] in the present paper
a new approach to geometrical quantization procedure is introduced.
As a first result we verify that the correspondence between
"classical" Poisson bracket on based symplectic manifold and
"quantum" Poisson bracket on the moduli space takes place.
\endabstract

   \endtopmatter

\head \S 0. Introduction \endhead                                               
$$
$$

In recent papers [2] and [8] it was represented construction of
moduli space of Bohr - Sommerfeld Lagrangian subcycles of fixed
topological type over a symplectic manifold. Starting with a simple
description of the tangent bundle to the moduli space authors propose
an additional superstructure --- fibered over the first moduli space
the second one of half - weighted Bohr - Sommerfeld Lagrangian subcycles 
of a fixed volume. This moduli space is endowed with structure of an infinite
dimensional kahler manifold. It's remarkable that this infinite dimensional manifold is constructed almost canonically over a fixed compact symplectic
manifold --- one has to add just a fixed metaplectic structure (see [2], [8]).

On the other hand one has an interesting approach to the description of
the quantum mechanics (see [1]). Instead of a Hilbert space usual in
the known descriptions one can use the projectivization of this space
as the phase space --- and this projectivization is an infinite dimensional
kahler manifold, endowed as a kahler manifold with an integrable complex structure, a symplectic 2- form and the corresponding Riemannian metric.
Then (see [1]) the Schroedinger equation describing the dynamics is replaced
by "classical" Hamilton equation using the symplectic structure. The Riemannian metric on the same projective space reflects the probabilistic aspects 
of the model. 

So dealing with these two constructions it is a natural idea
in the geometric quantization programme. Namely for a classical mechanical
system represented by the corresponding compact symplectic manifold $(M, \om)$
with symplectic form $\om$ one can consider an infinite dimensional
kahler manifold defined as the moduli space of half - weighted Bohr - 
Sommerfeld Lagrangian subcycles of a fixed volume. This infinite
dimensional kahler manifold can be regarded as an analogy of the projectivized
Hilbert space corresponds to an appropriate quantization of the system.
In this setup there is a natural correspondence between the space of smooth
functions on $(M, \om)$ and a special subspace in the space of smooth functions
over the moduli space. For a compact symplectic manifold $(M, \om)$
and the corresponding infinite dimensional kahler manifold $\sB^{hw,1}_{BS}$
one has an inclusion
$$
F: C^{\infty}(M \to \R) \to C^{\infty} (\sB^{hw, 1}_{BS}),
\tag 0.1
$$
explicitly described below in the formula (3.3). So the natural question arises:
what is the relationship between the Poisson brackets on the original finite dimensional symplectic manifold and on the induced infinite dimensional
kahler manifold which admits canonical symplectic structure?

The main result of the present paper is the following (Proposition 3.2):
inclusion $F$ preserves the structure of Lee algebra on $C^{\infty} (M \to \R)$,
and the Poisson bracket induced by $\om$ maps to the Poisson bracket induced
by the canonical symplectic form $\Om$ up to a constant multiple (at the end of section 4 we discuss what is the constant).

In other words if the infinite dimensional kahler manifold consists of
half - weighted Bohr- Sommerfeld Lagrangian subcycles of fixed
volume is a good candidate on the role of  quantum phase space
of the quantized system then the correspondence between the Poisson brackets
ensures us that the Dirac condition holds in this quantization procedure (about
the Dirac condition, prequantization, quantization etc. see [3], [4], [7]).

\head  1. Moduli spaces of half- weighted Bohr - Sommerfeld Lagrangian subcycles
\endhead

Let us recall briefly the main constructions and formulae from [2] and [8]
which we use in what follows.

Let $(M, \om)$ be a compact symplectic manifold of real dimension $2 n$
endowed with symplectic form $\om$. For an appropriate real smooth n- dimensional manifold $S$ one considers the space of all smooth inclusions
$$
\phi: S \to M
\tag 1.1
$$
such that $\phi^* \om$ is identically zero on $S$. It's easy to see
that this space just corresponds to the space of all smooth Lagrangian cycles
on $M$, homologically equivalent to $S$. Let moreover the pair $(M, \om)$
satisfies the integer Dirac condition
$$
[\om] \in H^2(X, \Z) \subset H^2(X, \R).
$$ 
Then one has on $(M, \om)$ a prequantization quadruple (see [7]) consists of
the manifold $M$,  the symplectic form $\om$ and additionally the complex line bundle $L \to M$ uniquely defined by the topological condition
$$
c_1(L) = [\om]
$$
and a hermitian connection $a \in \sA_h(L)$ which is a solution of the natural equation
$$
F_a = 2 \pi \imath \om.
\tag 1.2
$$
A Lagrangian inclusion $\phi$ (1.1) is Bohr - Sommerfeld if the following condition is satisfied: for the pair $(\phi^* L, \phi^* a)$ over $S$
where $\phi^* L$ is topologically trivial over $S$ and $\phi^* a$ is a flat connection (see (1.2))  the last one is gauge equivalent to the ordinary $d$.

It's easy to see that if $S$ is simply connected then every Lagrangian inclusion is Bohr - Sommerfeld.

After factorization of the space of such inclusion by the $Diff S$- action
one gets a space which was called the moduli space of Bohr - Sommerfeld Lagrangian subcycles of fixed topological type. In the pioneer papers [2] and
[8] the tangent space to the moduli space in a smooth point $(S, \phi)$ is represented by the space of exact 1-forms on $S$. For simplicity let us
fix an appropriate $S$ and denote, following the original papers, as 
$\sB_{BS}$ the moduli space. But the space $\sB_{BS}$ doesn't admit any
natural symplectic structure. 

So the next step in the construction to consider the space of Planckian cycles
$\Cal P_S$ fibered over $\sB_{BS}$ with $U(1)$- fibers. Take now the space 
of all half- weighted Planckian cycles $\Cal P^{hw}_{BS}$ consists of
pairs
$$
\Cal P^{hw}_{BS} = \{ (\rho, \theta) \},
$$
where $\rho$ is a Planckian cycle over $M$ and $\theta$ is a half - form over $\rho$ (see [2], [8]). The point is that the last space admits a natural symplectic structure and a natural compatible integrable complex structure
so it is a kahler manifold ([2], [8]). Further, on the kahler manifold
$\Cal P^{hw}_{BS}$ one has the natural $U(1)$- action, preserving both the structures. And we have as a result of kahler reduction (in dependence with
choosing fixed volume) the moduli space $\sB^{hw, r}_{BS}$ of half- weighted Bohr - Sommerfeld Lagrangian cycles with fixed volume $r \in \R$. In what follows we work with cycles of volume 1 just for simplicity but evidently all the results hold for any $r$.

For the computations below we need to recall the local description 
of $\sB^{hw,1}_{BS}$ and the explicit formula for $\Om$.  Points of
$\sB^{hw, 1}_{BS}$ are represented by pairs $(\phi, \theta)$ where $\phi$
is a Bohr - Sommerfeld Lagrangian inclusion of $S$ into $M$ and $\theta$
is a half - form on $\phi S \subset M$ such that
$$
\int_{\phi S} \theta^2 = 1.
$$
The tangent space $T_{(\phi_0, \theta_0)} \sB^{hw, 1}_{BS}$ is represented by pairs $(f, \theta)$ where $f$ - is a smooth function on $\phi_0 S$ and $\theta$ is a half- form on the same sub manifold such that
$$
\int_{\phi_0 S} f \theta_0^2 = 0
\tag 1.4
$$
and
$$
\int_{\phi_0 S} \theta_0 \theta = 0.
\tag 1.4'
$$
The symplectic form $\Om$ at the point $(\phi_0, \theta_0)$ has the form
$$
\Om(v_1, v_2) = \int_{\phi_0 S} (f_1 \theta_2 - f_2 \theta_1) \theta_0
\tag 1.5
$$
for every two tangent vectors
$$
v_1 = (f_1, \theta_1), \quad \quad v_2 = (f_2, \theta_2).
$$
It's not hard to see that this 2- form is nondegenereted everywhere
and is closed. All details can be found in [2] and [8].

\head  2. The background: geometrical formulation of quantum mechanics
\endhead

So we have seen that for every compact symplectic manifold $(M, \om)$ there exists a set of infinite dimensional kahler manifolds in dependence of
the homological classes. It looks very interesting if we turn to the framework
of "geometrical formulation of quantum mechanics". Since the author learned
this subject from the article [1] it is recommended as a reference together with [6]. 

Let $\Cal H$  is a Hilbert space (here and in what follows in the present section we will use the notations of [1], [6]) corresponding to a quantum mechanical system. Consider $\Cal H$ as a real vector space endowed with
complex structure $J$. Then the hermitian inner product can be decomposed 
into the sum of real and imaginary parts
$$
<\Phi, \Psi> = \frac{1}{2 \hbar} G (\Phi, \Psi) + \frac{\imath}{2 \hbar}
\Om(\Phi, \Psi).
\tag 2.1
$$
Thus on the real vector space $\Cal H$ we have the corresponding triple $(G, J, \Om)$, consists of positive defined real inner product, complex structure and the corresponding symplectic form.

On the Hilbert vector space every observable (being represented by a self adjoint operator) can be regarded as a vector field on $\Cal H$: really in each point of $\Cal H$ (more precisely  at each vector) we have another vector - the volume of this operator. The Schroedinger equation reads as
$$
\Psi' = - \frac{1}{\hbar} J \hat H \Psi,
\tag 2.2
$$
so it's convenient to use the following notation
$$
Y_{\hat F} (\Psi) = - \frac{1}{\hbar} J \hat F \Psi
\tag 2.3
$$
for the vector field corresponding to an observable $\hat F$. For the same
observable $\hat F$ let us consider the expectation value function
$$
F: \Cal H \to \R, \quad \quad F(\Psi) = <\Psi;, \hat F \Psi> =
\frac{1}{2 \hbar} G(\Psi; \hat F \Psi).
$$
It's easy to see ((2.5) in [1]) that
$$
i_{Y_{\hat F}} \Om = d F,
\tag 2.4
$$
so the Schroedinger equation reads as the Hamilton equation. Moreover, for two
observables $\hat F, \hat K$ the Poisson bracket of the functions
$F$ and $K$ is represented as the expectation value function of the self - adjoint operator equals to a multiple of their commutator (see (2.6), [1]).

The relationships of these types use to be  applicable after we turn from the
real kahler space to the projectivization of $\Cal H$, which was denoted as $\Cal P$ in [1]. Section B of the paper [1] is dedicated to a construction of
such projective space using the language of constraints. From the mathematical point of view it is equivalent to construct the projective space using 
the kahler reduction under the natural $U(1)$- action with moment map
$$
S: \Cal H \to \R, \quad \quad S(\Psi) = <\Psi; \Psi>
\tag 2.5
$$
for any level, for example for  level $S(\Psi) = 1$. Under the projective space $\Cal P$, realized by this procedure, one has the corresponding kahler structure
consists of Riemannian metric $g$, integrable complex structure $I$ and
symplectic 2-form $\Om$. Instead of observables $\hat F$ one could consider
the induced expectation value functions $F$ which are invariant under the $(U(1)$- action being restricted on the hypersurface $S(\Psi) = 1$ thus these function are correctly defined on $\Cal P$. In section C of the paper [1]
one shows that the correspondence between usual observables on $\Cal H$ and
 special functions on $\Cal P$ is one - to - one, so one can reconstruct 
the original self - adjoint operator from the corresponding expectation function
on $\Cal P$. On the other hand one has to specify the subspace of observable functions over $\Cal P$, induced by self - adjoint operators. The answer is:
 smooth function
$$
f: \Cal P \to \R
$$
 is induced by a self - adjoint operator if and only if its Hamiltonian vector
field $H_f$ over $\Cal P$ is a Killing vector field for the Riemannian metric $g$ (Corollary 1, Theorem 2.1 in [1]0.

Before  forgetting about the original Hilbert space it's necessary to translate
the notions "eigenvector" and "eigenvalue" to the projective language. It's not hard to see (subsection 3, section C of [1]) that

--- for observable $\hat F$ with induced function $f: \Cal P \to \R$ every eigenvector after the projectivization uses to be a critical point of the function $f$;

--- and the corresponding critical value equals to the original eigenvalue.

Now we are ready to formulate the postulates of the quantum mechanics in these projective terms following the authors of [1] (see section D [1]).
Quantum phase space is represented by an appropriate projective space $\Cal P$ with is a kahler manifold (finite dimensional or infinite dimensional); the space of observable is the set of real smooth functions on $\Cal P$ such that
their Hamiltonian vector fields are Killing vector fields with respect to the
Riemannian metric; the dynamics is described by the Hamilton equation;
all probabilistic aspects, state reductions and so on are described in terms
of geodesic distances (so based on the Riemannian structure on $\Cal P$). One has the state reductions in both cases (of discrete spectrum and of non isolated
critical points).

Roughly speaking the authors show that the difference between the classical mechanics and the quantum one is in the presence of an appropriate Riemannian metric. In both cases one has symplectic structures, reflecting
the dynamics, but additionally in the quantum case one has a Riemannian metric,
reflects the probabilistic aspects. At the end of the paper the author
propose the following question: is there exist a quantization procedure
of classical mechanical systems  which gives directly from a given classical system an appropriate infinite dimensional kahler manifold $\Cal P$ with Riemannian metric $g$ and symplectic 2- form $\Om$ and doesn't use known approaches, deriving form the system Hilbert vector spaces and doesn't refer
to such Hilbert spaces?

\head  3. Preferred functions on $\Cal B^{hw,1}_S$
\endhead

In this section we construct the inclusion $F$ (0.1) and consider some geometrical objects on $\sB^{hw, 1}_S$ induced by a smooth real function $f$ 
defined on the based symplectic manifold $(M, \om)$. 

Let $f \in C^{\infty}(M \to \R)$ is a smooth function. Then the differential
of this function $df$ being restricted on a Bohr - Sommerfeld cycle $\phi S
\subset M$ gives us a tangent vector to $\sB^{hw, 1}_S$ in point $(\phi, \theta)$ for any $\theta$. This gives a vector field on $\sB^{hw, 1}_S$
denoted as $A_f$. This vector field doesn't depend on the second "coordinate"
on $\sB^{hw, 1}_S$, in other words it is constant along fibers of
$$
\sB^{hw, 1}_S \to \sB_S.
$$
We have the following simple
\proclaim{Proposition 3.1} The set of singular points of $A_f$ consists of
such Bohr - Sommerfeld cycles $\phi S$ that the function $f$ is constant being restricted on $\phi S$. 
\endproclaim

Together with the vector field for a function $f$ one has a natural induced 1- form  on $\sB^{hw, 1}_S$. At a point $(\phi_0, \theta_0)$ this form denoted as
$B_f$ reads as
$$
B^f_{(\phi_0, \theta_0)}(f_1, \theta_1) = \int_{\phi S} f \theta_0 \theta_1.
\tag 3.1
$$
Since the symplectic form $\Om$ is described by the formula (1.5) one gets 
by direct computation that the vector field $A^f$ and the 1- form $B_f$
are related as follows
$$
B_f = \Om^{-1}(A^f).
\tag 3.2
$$

Now we want to define the inclusion $F$. For each smooth function
$f$ on $M$ let us define the following function $F_f$, naturally induced on
$\sB^{hw, 1}_S$:
$$
F_f (\phi, \theta) = \int_{\phi S} f|_{\phi S} \theta^2.
\tag 3.3
$$

First of all let us remark that the images of constant functions are
constant functions on $\sB^{hw, 1}_S$ (compare (3.3) and (1.3)).

Consider the image of the smooth function space on $M$ as a subspace of
the smooth function space on $\sB^{hw,1}_S$:
$$
\Cal N = Im_F (C^{\infty} (M \to \R)), \quad \Cal N \subset C^{\infty}(\sB^{hw, 1}_S \to \R).
$$
It's easy to see that $\Cal N$ is a linear subspace. But if we consider 
the space $C^{\infty}(\sB^{hw, 1}_S \to \R)$ as the algebra with pointwise
multiplication then the subspace $\Cal N$ were not a subalgebra. Really if
$F_f$ and $F_g$ are two induced functions form $\Cal N$ then 
the product $F_f \cdot F_g$ isn't a priori induced by a real function and
doesn't lie in $\Cal N$. But this shortage is compensated by the fact that
this $\Cal N$ is a Lie subalgebra with respect to the Poisson bracket induced by
the symplectic structure $\Om$. We have the following 
\proclaim{Proposition 3.2} The identity
$$
\{ F_f; F_g \}_{\Om} = 2 F_{\{f; g\}_{\om}}
$$ holds.
\endproclaim

Proposition 3.2 together with definition (3.3) represent an answer
(or just a part of an answer) on the question arises in Introduction. Before
we will prove the statement let us input  few remarks.

The definition (3.3) of the correspondence $f \mapsto F_f$ can be easily modified scaling by an appropriate constant
$$
F_f (\phi_0, \theta_0) = \int_{\phi_0 S} \tau f|_{\phi_0 S} \theta_0^2.
$$
This modification changes the identity (3.5) as follows
$$
\{ F_f; F_g \}_{\Om} = 2 \tau^2 F_{\{f; g \}_{\om}}.
$$
So one can rearrange the identity (3.5) such that any physical constant
(Planck constant etc.) will be consistent. Thus the Dirac condition 
is realizable in the framework.

On the other hand, it has been remarked that the pointwise multiplication
doesn't preserve the subspace $\Cal N$. This uninvariance gives us an interesting effect. Let us suppose that the given classical mechanical system, represented by symplectic manifold $(M, \om)$, is completely integrable.
Thus we have a number of the integrals $f_1, ..., f_n$, which
 commute each with others. This set generates a subalgebra $\Cal V \subset C^{\infty}(M \to \R)$ with usual pointwise multiplication, and this $\Cal V$
is the maximal commutative Lie subalgebra  in $C^{\infty}(M \to \R)$, viewing as a Lie algebra. In other words the maximal commutative Lie subalgebra is 
finitely generated as usual algebra: for any element $h$ from this maximal Lie subalgebra there are exist  a $n$- tuple $(r_1, ..., r_n), r_i \in \Bbb N$ such that 
$$ 
h = f_1^{r_1} \cdot ... \cdot f_n^{r_n}.
$$
Let us turn now to the infinite dimensional moduli space $\sB^{hw, 1}_S$. 
For the completely integrable situation one has the following infinite set of
commuting function
$$
\{ F_{f_1^{r_1} \cdot ... \cdot f_n^{r_n}} \},
$$
but for this maximal commutative Lie subalgebra the same isn't true ---
we have at least $\Bbb N^n$  generators because $F_{f_1 \cdot f_2} \ne
F_{f_1} \cdot F_{f_2}$. So one could define the dimension of symplectic manifold $\sB^{hw, 1}_S$ using this arguments.

\head 4. Computations
\endhead

Here we prove Proposition 3.2. Let $f$ is a smooth function on $M$.
If we take the corresponding $F_f \in C^{\infty}(\sB^{hw, 1}_S)$
then it has the differential equals to
$$
dF_f (\phi_0, \theta_0) (f_1, \theta_1) = \int_{\phi_0 S} 2 f \theta_1 \theta_0
+ \int_{\phi_0 S} df_1 (\om^{-1}(df)|_{\phi_0 S}) \theta_0^2,
\tag 4.1
$$
so the first summand in (4.1) is constant under the varying of the first "variable" and the second summand is constant under the varying of the second one. One could recognize the form $B_f$ defined by (3.1) as the first summand in (4.1). As we have seen above (see (3.2)) the vector field $A^f$ is symplectically dual to the form. Thus the Hamiltonian vector field of the function $F_f$ has the form
$$
H_{F_f} = \Om^{-1} (d F_f) = 2 A^f + C^f,
$$
where $C^f$ is the vector field, symplectically dual to the 1- form
$$
\int_{\phi_0 S} df_1 (\om^{-1}(df)|_{\phi_0 S}) \theta^2_0.    \tag 4.2
$$
we need not to compute this $C^f$ explicitly because of the following argument.
Vector field $A^f$ is constant on the half form coordinates while the symplectically dual 1- form $B_f$ has as the kernel all vector field of the shape $(f_1(\phi, \theta), 0)$. Vice versa for the vector field  $C^f$ and for the symplectic dual 1- form one has the same properties if one changes coordinates "functions --- half forms". From the formula (1.5) we see that the symplectic form $\Om$ separates the coordinates $f$ and $\theta$. So one can see that
$$
\aligned
\{ F_f, F_g \}_{\Om} = \Om(H_{F_f}; H_{F_g}) = \\
\Om (2 A^f + C^f; 2 A^g + C^g) = 2 \Om (A^f; C^g) - 2 \Om (A^g; C^f) \\
\endaligned
\tag 4.3
$$
since $\Om(A^f; A^g) = \Om (C^f; C^g) = 0$ for every pair $f, g$. Further,
as a continuation of (4.3) one gets
$$
\{ F_f; F_g \}_{\Om} = 2 B_f(C^g) - 2 B_g (C^f) = 2 C^*_f (A^g) - 2 C^*_g (A^f),
\tag 4.4
$$
where we denote  the form (4.2) as $C^*_f$. Thus we need not to compute the explicit expression for $C^f$. Really substituting the explicit expressions for the vector fields and 1- forms we get
$$
\{ F_f; F_g \}_{\Om} = 2 \int_{\phi_0 S} dg |_{\phi_0 S} (\om^{-1}(df)|_{\phi_0 S}) \theta_0^2 - 2 \int_{\phi_0 S} df|_{\phi_0 S} (\om^{-1}(dg)|_{\phi_0 S})
\theta_0^2.
\tag 4.5
$$
Let us take the integrand from (4.5)
$$
2 (dg|_{\phi_0 S}(\om^{-1}(df)|_{\phi_0 S}) - df|_{\phi_0 S} (\om^{-1}(df)|_{\phi_0 S})).
\tag 4.6
$$
It is a function on $\phi_0 S$. Let us show that this function coincides
with the restriction on $\phi_0 S$ of the Poisson bracket $\{f; g \}_{\om}$, multiplied by 2. The rest of this section is dedicated to the proof of this coincidence.

First of all it's easy to see that the following identity takes place
$$
2 \{ f; g \}_{\om} = 2 df (\om^{-1}(dg)) = - 2 dg (\om^{-1}(df)) =
(df (\om^{-1}(dg)) - dg (\om^{-1}(df)).
\tag 4.7
$$
But we work now near a Lagrangian submanifold --- we are interesting in the expression (4.7) only in a neighborhood of our $\phi_0 S$. For simplicity let us choose any compatible almost complex structure $J$ on $M$, getting the corresponding hermitian triple $(g, J, \om)$ on the based symplectic manifold.
The corresponding Riemannian metric $g$ splits near our Lagrangian submanifold
$\phi_0 S$ which means that
$$
v \in T_m \phi_0 S \subset T_m M \implies G(v, Jv) = 0.
\tag 4.8
$$
In the expression (4.6) the restrictions on $\phi_0 S$ of 1- forms $df, dg$
and vector fields $\om^{-1}(dg), \om^{-1}(df)$ take place. So let us decompose
every ingredient in  formula (4.7) into horizontal and vertical parts with respect to the tangents subspaces to $\phi_0 S$ and their orthogonal complements
with respect to the Riemannian metric $g$. For the first summand in (4.7) one gets
$$
\aligned
(df_{vert}((Jg^{-1}(dg))_{vert}) + df_{vert}((Jg^{-1}(dg))_{hor}) + \\
df_{hor}((Jg^{-1}(dg))_{vert}) + df_{hor}((Jg^{-1}(dg))_{hor})).
\\
\endaligned
\tag 4.9
$$
In the last expression (4.9) one has only two nontrivial summands due to the orthogonality (4.8) --- namely the first and the forth ("vert - vert" and
"hor - hor"). Analogously for the second summand in (4.7) one has
$$
- (dg_{vert}((Jg^{-1}(df))_{vert}) + d g_{hor}((Jg^{-1}(dg))_{hor})),
\tag 4.10
$$
so we have again only two nontrivial terms. But from the compatibility condition for the Riemannian metric and the almost complex structure one has
$$
df_{vert}((J g^{-1}(dg))_{vert}) = - dg_{hor}((J g^{-1}(df))_{hor})
\tag 4.11
$$
and as well
$$
df_{hor}((Jg^{-1}(dg))_{hor}) = - dg_{vert}((Jg^{-1}(df))_{vert}).
\tag 4.12
$$

Thus we can rewrite the expression for the Poisson bracket $\{f; g \}_{\om}$ restricted on the Lagrangian submanifold $\phi_0 S$  using only horizontal components. This gives us
$$
\{ f; g \}_{\om}|_{\phi_0 S} = df_{hor}((\om^{-1}(dg))_{hor}) - dg_{hor}((\om^{-1}(df))_{hor}).
\tag 4.13
$$
Now it remains  to recall that the $hor$ - components  just correspond to the restrictions of vector fields and 1- forms to the Lagrangian submanifold $\phi_0 S$ thus we have the following identity
$$
\{f; g \}_{\om} |_{\phi_0 S} = (df|_{\phi_0 S}(\om^{-1}(dg)|_{\phi_0 S})
- dg|_{\phi_0 S}(\om^{-1}(df)|_{\phi_0 S}).
\tag 4.14
$$
Comparing (4.6) and (4.14) one gets the statement of Proposition 3.2.

\head Final remarks
\endhead

This paper contains  first results about the moduli spaces of half - weighted Lagrangian subcycles of fixed volume. We got an interesting result about
the relationship between the Poisson structures on the based manifold and on the moduli space, nothing else. But in the framework on the geometrical quantization programme this result looks like a hint that the approach mentioned in the introduction to this paper could be exploited. Really, as usual one understands
the geometric quantization of a classical mechanical system as a procedure
gives an appropriate Hilbert space together with a correspondence "classical observables --- quantum observables" such that the irreducibility and the Dirac conditions would be
satisfied (see [3], [4], [7] etc.). Leaving the question of the possibility 
to construct such correspondence (we mean the celebrated van Hove theorem, see f.e.
section 5.2 in [3])let us turn to the geometrical formulation of quantum mechanics. As we have seen in section 2 (following A. Ashtekar and T. Schilling)
one could try to find an appropriate kahler manifold instead of a Hilbert space
directly. It has to have some very special properties (maximal symmetries and so on) but as well it has to be defined directly from the based symplectic manifold. At the present time we have a candidate --- the moduli space of half- weighted Bohr - Sommerfeld Lagrangian subcycles constructed by A.Tyurin and
A.Gorodentsev. Their moduli space has to be studied in the  following directions:

1. irreducibility of the moduli space (for the irreducibility condition);

2. maximal symmetries (to be really a geometrical model of quantum mechanics);

3.  the functions $F_f$ from $\Cal N$ (to be really quantum observables).

Preliminary results, which we got during the work on this paper (but
were not included here because of the negativity), show that at least for the third question the answer has to be negative. The point is that for any $F_f$
the critical points form a continuous set, and there are exist such functions 
from $\Cal N$ which have only one critical value. But the situation is very rich
to avoid a lot of difficulties

--- we can study some special submanifolds of $\sB^{hw,1}_S$ to get
a more appropriate kahler manifold; this way can give us
right constructions to move in the 1st and the 2nd directions;

--- we can modify the definition of the correspondence $F$ (0.1)
to move in the 3d direction;

--- at the end, we have seen that our result (Proposition 3.2)
depends on the symplectic structure on $\sB^{hw,1}_S$ and doesn't depend 
on Riemannian metric or integrable complex structure. So we can perturb the original complex structure together with the Riemannian metric (constructed in [2], [8]) getting an almost complex symplectic manifold. May be this almost complex manifold gives the right construction of geometric quantization. 
If this approach will be succesfull  we would relax the celebrated Penrose slogan: "The Nature is complex". We will add:"or almost complex".

At the end I would like to express my gratitude to Max - Planck Institute fur Matematik (Bonn) for hospitality. As well I have to thank A.Gorodentsev, P. Pyatov and P.Saponov for valuable discussions and remarks.

\enddocument